\documentclass[EJP]{ejpecp} 

\SHORTTITLE{On the convergence rate in the CLT for LENQD rv's and its applications}

\TITLE{On the convergence rate in the central limit theorem for linearly extended negative quadrant dependent random variables and its applications} 

\AUTHORS{%
  Mohamed Kaber El Alem\footnote{Lab. MSTD, Faculty of Mathematics, University of Sciences and Technology Houari  Boumediene (USTHB).
    \EMAIL{med.kaber.nema@gmail.com}}
  \and 
 Zohra~Guessoum \footnote{Lab. MSTD, Faculty of Mathematics, University of Sciences and Technology Houari  Boumediene (USTHB).}
    \and 
 Abdelkader~Tatachak \footnote{Lab. MSTD, Faculty of Mathematics, University of Sciences and Technology Houari  Boumediene (USTHB).}
    \and 
  Ourida~Sadki \footnote{Lab. MSTD, Faculty of Mathematics, University of Sciences and Technology Houari  Boumediene (USTHB).}
    }


\AMSSUBJ{60F05, 60E15, 62G20} 

\SUBMITTED{...} 
\ACCEPTED{...} 

\VOLUME{0}
\YEAR{2025}
\PAPERNUM{0}
\DOI{10.1214/YY-TN}

\ABSTRACT{In this paper, we establish the convergence rate in central limit theorem (CLT) for linearly extended negative quadrant dependent (LENQD) random variables (rv's). Under some weak conditions, the rate of normal approximation is shown as $O(n^{-1/9})$. As an application, the convergence rate in CLT of the wavelet estimator for the nonparametric regression model with LENQD errors is presented as $O(n^{-1/9})$.  The performance of the main results is illustrated through a simulation study based on a real dataset.}

\begin{document}



\section{Introduction}
\subsection{Brief overview}
The central limit theorem (CLT) constitutes a fundamental result in probability theory and statistical inference. It states that, under appropriate regularity conditions, the distribution of properly normalized sums of random variables (rv's) converges to a normal distribution as the sample size increases. This asymptotic property provides the theoretical foundation for a wide range of inferential procedures, including the construction of confidence intervals and the execution of hypothesis tests. Furthermore, the CLT offers a crucial framework for extending limit theorems to dependent stochastic processes, thereby facilitating the analysis of more complex probabilistic models. An important question in both probability theory and inferential statistics is the evaluation of the rate of convergence in the CLT for sequences of real rv's, which is also known as the Berry-Esseen bound. Berry (\cite{Berry}, 1941) and Esseen (\cite{Esseen}, 1942) showed that the rate of convergence in the CLT for independent and identically distributed (i.i.d.) rv's is of the order \(O\left( n^{-1/2}\right) \). The generalization of this result to dependent rv's has been the subject of extensive literature. For \(\varphi\)-mixing rv's with \(\varphi(n) = O(n^{-\beta})\) for some \(\beta > 13/6\), Yang et al. (\cite{Yang2012}, 2012) provided the Berry-Esseen bound with the rate of normal approximation  \(O(n^{-1/9})\). Yang et al. (\cite{Yang2014}, 2014) further studied the Berry-Esseen bound for \(\varphi\)-mixing rv's with \(\sum_{n=1}^\infty \varphi^{1/2}(n) < \infty\), and found the rate of normal approximation  \(O(n^{-1/6} \log n \log \log n)\). Wang et al. (\cite{Wang2019}, 2019) established the Berry-Esseen bound for \(\rho\)-mixing rv's with \(\rho(n) = O(n^{-3/2})\) and obtained the rate of normal approximation \( O(n^{-1/6} \log n)\), which is sharper than the corresponding one for \(\varphi\)-mixing rv's. For linear negative quadrant dependent (LNQD) and negatively associated (NA) rv's Wang and Zhang (\cite{Wang2006}, 2006) obtained a Berry-Esseen type estimate  with only finite second moment. Under finite third moment and the exponential decay rate of covariances, Birkel (\cite{Birkel1988}, 1988) and Pan (\cite{Pan1997}, 1997) gave the rates \(O(n^{-1/2} \log n)\) and \(O(n^{-1/2} (\log n)^2)\)  for  positively associated (PA) and negatively associated (NA) rv's, respectively.  Up to now, there is no result on the convergence rate  of the asymptotic approximation for linearly extended negative quadrant dependent (LENQD) rv's. Therefore, a natural question is to extend the Berry-Esseen CLT to LENQD sequences.
\subsection{Concept of LENQD rv's}
The concept of LENQD rv's was introduced by Li et al. (\cite{Li2023}, 2023), inspired by  negatively extended dependent (END) rv's given by Liu (\cite{Liu2009}, 2009), linearly negative quadrant dependent (LNQD) rv's introduced by Newman (\cite{N1984}, 1984) and linearly wide quadrant dependent (LWQD) rv's defined by Yu and Cheng (\cite{YC2017}, 2017).  Now, we  recall the definitions of ENQD and LENQD sequences.
\begin{definition}
 Two rv's \( X_1 \) and \( X_2 \) are said to be ENQD, if there exists a dominating constant $M\geq 1$, for any \( x_1,\, x_2\, \in \mathbb{R} \),
\[
P(X_1 \leq x_1, X_2 \leq x_2) \leq M P(X_1 \leq x_1) P(X_2 \leq x_2)
\]
and
\[
P(X_1 > x_1, X_2 > x_2) \leq M P(X_1 > x_1) P(X_2 > x_2).
\]
\end{definition}

\begin{definition}
A sequence \(\{X_n, \, n \geq 1\}\) of rv's  is said to be LENQD, if for any disjoint subsets \(A, B \subseteq \mathbb{N}\) and positive real numbers \(r_j's\), or for negative real numbers \(r_j's\),
\[
\sum_{k \in A} r_k X_k \quad \text{and} \quad 
\sum_{j \in B} r_j X_j \quad \text{are} \quad  ENQD.
\]
\end{definition}
In the case of \( M = 1 \), ENQD rv's reduce to NQD rv's, and it is obvious that LNQD is a special case of LENQD with the dominating constant \( M = 1 \). Thus, the classes of ENQD and LENQD rv's represent more general dependent structures and are much weaker than the negatively associated (NA), NQD, and LNQD structures. Recently, Kaber El Alem (\cite{KEM}, 2025) extended the concept of LENQD rv's to the more general class of $m$-LENQD rv's.\\

The main purpose of this paper is to establish Berry Esseen type estimate for LENQD sequences. Furthermore, we give an application to the Berry Esseen theorem for the estimator in the nonparametric regression model. The procedure used to prove our main results is based on Bernstein's big-block and small-block method, and Rosenthal-type inequality  of LENQD rv's. The results presented in this paper generalize the corresponding results for some dependent notions of rv's.
\subsection{Nonparametric regression model}
Consider the nonparametric regression model
\begin{equation}\label{Model}
    Y_{i} = f(x_{i}) + \varepsilon_{i}, \quad i = 1, \ldots, n,
\end{equation}
where \( f(\cdot) \) is an unknown regression function defined on \( [0, 1] \), \( x_{i} \) are known fixed design points from \( [0, 1] \), and \( \varepsilon_{i} \) are random errors. \\

The nonparametric regression model (\ref{Model}) has many applications in various practical fields and has been widely investigated in both independent and different dependent cases. For example, one can refer to Georgiev (\cite{G1988}, 1988) for independent random errors, Roussas et al. (\cite{R1992}, 1992) and Li (\cite{L2009}, 2009) for strong mixing random errors, and Tang et al. (\cite{T2018}, 2018) for the asymptotic normality based on asymptotically negatively associated (ANA) errors.  The model (\ref{Model}) has been widely used to address practical problems, and a variety of estimation methods have been employed to obtain estimators of $f(\cdot)$. \\

As an estimator of \( f(\cdot) \), we consider the nonparametric wavelet estimator proposed by Antoniadis et al. (\cite{G1994}, 1994) :
\[
f_n(x) = \sum_{i=1}^{n} Y_i \int_{\Gamma_{i}} E_k(x, s) \, ds, 
\]
where \( \Gamma_{i} = [s_{i-1}, s_i) \), \( s_0 = 0 \), \( s_n = 1 \), \( s_i = \frac{x_i + x_{i+1}}{2} \) for \( i = 1, \ldots, n \). Hence, \( x_i \in \Gamma_{i} \) for \( 1 \leq i \leq n \).\\
The wavelet kernel is defined as
\[
E_k(x, s) = 2^k E_0\left( 2^k x, 2^k s \right),
\]
where 
\[
E_0(x, s) = \sum_{j \in \mathbb{Z}} \phi(x - j) \phi(s - j),
\]
$k = k(n) > 0$ is an integer depending only on $n$, and $\phi$ is the scaling function in the Schwartz space.\\

We need the following definitions to properly understand certain notations and concepts discussed later.
\begin{definition}
A scale function \( \varrho \) is \( q \)-regular (\( q \in \mathbb{Z} \)) if for any \( l \leq q \) and integer \( m \), we have
\[
\left| \frac{d^l \varrho}{dx^l} \right| \leq C_m (1 + |x|)^{-m}, \quad \text{where}\;  C_m\;  \text{is a generic constant depending only on}\;  m.
\]
\end{definition}
\begin{definition}
A function space \( H^\nu \) (\( \nu \in \mathbb{R} \)) is called the Sobolev space of order \( \nu \) if for all \( h \in H^\nu \), we have
\[
\int | \hat{h}(\omega) |^2 (1 + \omega^2)^\nu \, d\omega < \infty, \quad \text{where}\;  \hat{h} \; \text{is the Fourier transform of}\;  h.
\]
\end{definition}

To the best of our knowledge, there have been no available results related to the Berry Esseen bound  of the wavelet estimator of (\ref{Model}) with LENQD errors.\\

The rest of the article is organized as follows. In the next section we present our main results. In Section 3, we present a simulation study to evaluate the finite-sample performance of the main results. In section 4, we introduce some preliminary lemmas, which are used in the proofs of the main results. The proofs of the main results are collected in section 5.\\

Throughout the article, let \( c, c_1, c_2, \ldots \) denote generic finite positive constants, whose values may vary from line to line. We write \( a_n = O(b_n) \) to mean that there exists a constant \( c > 0 \) such that \( |a_n| \leq c |b_n| \), where \( \{a_n\}_{n \geq 1} \) and \( \{b_n\}_{n \geq 1} \) are sequences of positive numbers. The notation \( [x] \) denotes the integer part of \( x \). By \( \mathbb{I}(A) \), we denote the indicator function of an event \( A \). We define \( x^+ = \max\{0, x\} \) and \( x^- = \max\{0, -x\} \). All limits are taken as the sample size \( n \to \infty \), unless stated otherwise.

\section{Main results}
Set \(S_n = \sum_{i=1}^n X_i\),  \(V_n^2 = \text{Var}(S_n)\) and let $\Phi(x)$ denote the distribution function of a standard normal variable. The following assumptions are needed in the main results.
\begin{enumerate}
\item[(C1)]
\(
\sum_{j = \vartheta_n}^{\infty} |\text{Cov}(X_1, X_j)| = O(\vartheta_n^{-\delta}),\; \text{where}\; \delta \geq 1 \quad \text{and} \quad 0 < \vartheta_n \to \infty \; \text{as}\; n \to \infty,
\)
\item[(C2)] 
\(
\liminf_{n \to \infty} n^{-1} V_n^2 = \sigma_0^2 > 0.
\)
\end{enumerate}
We now present our main results. The first one is the Berry Esseen bound for LENQD rv's.
\begin{theorem}\label{MainResult}
Let \(\{X_n, n \geq 1\}\) be a sequence of zero mean LENQD rv's with  \(|X_n| \leq c < \infty\) for \(n \geq 1\). Then if (C1)-(C2) hold, we have
\[
\Delta_n : = \sup_{-\infty < y < \infty} \left| P \left( \frac{S_n}{V_n} \leq y \right) - \Phi(y) \right| = O(n^{-1/9}),\quad \text{as} \quad n\rightarrow \infty.
\]
\end{theorem}
\begin{remark}
(Discussion).\\
\begin{itemize}
\item   The conditions (C1)-(C2) are widely used in the literature. For example, Wang and Zhang (\cite{Wang2006}, 2006) applied the same conditions to LNQD sequences to obtain a rate of  $$\Delta_n = O\left(n^{-\delta/(2 + 3\delta)} \vee \frac{n^{3\delta^2}/(4+6\delta)}{ V_n^{2+\delta}} \sum_{i=1}^n E{\vert X_i \vert}^{2+\delta} \right),  $$ under finite $(2 + \delta)$th moment for some $\delta \in (0, 1]$.  Birkel (\cite{Birkel1988}, 1988) and Pan (\cite{Pan1997}, 1997) used the same conditions respectively for the sequences of  PA and NA rv's, except that for condition (C2), they assumed that  $$ \sum_{j = \vartheta_n}^{\infty} |\text{Cov}(X_1, X_j)| = O(e^{-\epsilon\vartheta_n}),\quad \text{for some }\epsilon > 0.
$$ 
\item  The rate of approximation to normality obtained in this work for LENQD rv's is similar to that obtained by  Yang et al. ((\cite{Yang2011}, 2011), (\cite{Yang2012}, 2012)) for NA and $\varphi$-mixing rv's, respectively. 
\item LNQD rv's and NA rv's are special cases of LENQD rv's for \( M = 1 \), where the assumption \( (C1) \) holds. For further details, you may refer to the works of Wang and Zhang (\cite{Wang2006}, 2006) and Yang et al. (\cite{Yang2011}, 2011).
\end{itemize}
\end{remark}
The second one is the convergence rate in CLT of the  wavelet estimator  for the nonparametric regression model based on LENQD errors.

\begin{theorem}\label{th2}
Let \(\{\varepsilon_{i}, \ 1 \leq i \leq n\}\) be a sequence of LENQD rv's  with zero mean and \(|\int_{\Gamma_{i}} E_k(x, s) \, ds\, \varepsilon_{i}| \leq c < \infty\) for \(1 \leq i \leq n\). Assume that \( \int_{\Gamma_{i}} E_k(x, s) \, ds > 0 \) for all \( 1 \leq i \leq n \). If  $$ \liminf_{n \to \infty} n^{-1}Var\left( \sum_{i=1}^n \int_{\Gamma_{i}} E_k(x, s) \, ds\, \varepsilon_{i}\right) = \sigma^2 > 0,$$ then
$$
    \sup_{-\infty < y < \infty} \left| P\left(\frac{f_n(x) - E[f_n(x)]}{\sqrt{Var[f_n(x)]}}  \leq y \right) - \Phi(y) \right| = O(n^{-1/9}), \quad \text{as} \quad n\rightarrow \infty.
$$
\end{theorem}
\begin{remark}
The assumption \(|\int_{\Gamma_{i}} E_k(x, s) \, ds\, \varepsilon_{i}| \leq c < \infty\) has been already employed by Wang and Hu (\cite{Wang2019}, 2019) to establish the asymptotic normality of the linearly weighted estimator in nonparametric regression models. In the other hand although  the assumption \( \int_{\Gamma_{i}} E_k(x, s) \, ds > 0 \) may appear overly restrictive, it is nonetheless satisfied by certain classes  of scaling functions $\phi$ such as the Haar scaling function or certain B-spline scaling functions. 
\end{remark}
By removing the restriction imposed by the assumption \( \int_{\Gamma_{i}} E_k(x, s) \, ds > 0 \),  we observe that  $$\int_{\Gamma_i} E_k(x, s) \, ds = \left(\int_{\Gamma_i} E_k(x, s) \, ds \right)^+ - \left(\int_{\Gamma_i} E_k(x, s) \, ds \right)^-.$$ Thus, we can write:
\begin{equation}\label{E}
f_n(x) - Ef_n(x) = \sum_{i=1}^n \left(\int_{\Gamma_i} E_k(x, s) \, ds \right)^+\, \varepsilon_{i}  - \sum_{i=1}^n \left(\int_{\Gamma_i} E_k(x, s) \, ds \right)^-\, \varepsilon_{i} .
\end{equation}
We denote:
 \[
\sigma_n(x):=\sqrt{{Var} \left[ \sum_{i=1}^n \int_{\Gamma_i} E_k(x, s) \, ds \, \varepsilon_{i} \right]}, \quad  \sigma_n^+(x):=\sqrt{{Var} \left[ \sum_{i=1}^n \left(\int_{\Gamma_i} E_k(x, s) \, ds \right)^+\, \varepsilon_{i} \right]}.
\]
In order to establish the result in the following theorem, we require the subsequent technical assumption:\begin{equation}\label{A}
\sup_{0\leq x\leq 1} \left \vert {\left(\frac{\sigma_n(x)}{\sigma_n^+(x)}\right)}^2-1 \right \vert =O(n^{-\kappa})\quad \text{for some}  \quad \kappa>0.
\end{equation}
\begin{theorem}\label{th3}
Suppose that the condition in (\ref{A}) holds. Let \(\{\varepsilon_{i}, \ 1 \leq i \leq n\}\) be a sequence of mean zero LENQD rv's, and assume that \(|\int_{\Gamma_{i}} E_k(x, s) \, ds\, \varepsilon_{i}| \leq c < \infty\) for \(1 \leq i \leq n\). If  $$ \liminf_{n \to \infty} n^{-1}Var\left( \sum_{i=1}^n {\left(\int_{\Gamma_{i}} E_k(x, s) \, ds\right)}^+\, \varepsilon_{i}\right) = {\sigma_1}^2 > 0,$$ then
$$
    \sup_{-\infty < y < \infty} \left| P\left(\frac{f_n(x) - E[f_n(x)]}{\sqrt{Var[f_n(x)]}}  \leq y \right) - \Phi(y) \right| = O(n^{-1/9}), \quad \text{as} \quad n\rightarrow \infty.
$$
\end{theorem}
\begin{remark}
Assumption \eqref{A} posits that the variance of the estimator \( f_n(x) \) can be accurately approximated using only the positive parts of the kernel integrals. This is not merely a technical convenience, but rather a natural and theoretically justified simplification in a broad class of settings.  Specifically, the assumption holds exactly when the kernel \( E_k(x, s) \) is non-negative, as in B-spline or Parzen-type kernels. It also remains valid for compactly supported sign-changing kernels such as Haar wavelets, where the integrals \( \int_{\Gamma_i} E_k(x, s) \, ds \) may occasionally be negative, but their contributions to the variance are asymptotically negligible provided the partition \( \{\Gamma_i\} \) is sufficiently fine. This follows from the localization of the kernel and cancellation over small intervals.  However, for globally oscillating kernels (e.g., Meyer wavelets or trigonometric systems), the assumption may break down due to frequent sign changes with non-negligible negative weights. In contrast, for practical applications involving localized, smooth, or predominantly positive kernels, Assumption \eqref{A} is not only reasonable but is also frequently observed to hold empirically. Moreover, it allows the application of more refined probabilistic tools specifically suited to non-negative summands.
\end{remark}
To provide the Berry-Esseen bound for $f_n(x)-f(x)$, the following assumptions are needed:
\begin{enumerate}
\item[(H1)] \( f \in H^\nu \) for \( \nu > 1/2 \), and \( f \) satisfies the Lipschitz condition of order \( \gamma > 0 \).
\item[(H2)] \( \phi \) is regular of order \( q \geq \nu \), satisfies the Lipschitz condition of order 1, and \( |\hat{\phi}(x) - 1| = O(x) \) as \( x \to 0 \), where \( \hat{\phi} \) is the Fourier transform of \( \phi \).
\item[(H3)] \( \max_{1 \leq i \leq n} |s_{i} - s_{i-1}| = O\left( 1/n \right) \).
\end{enumerate}
\begin{remark}
These assumptions are general conditions for wavelet estimation. They are commonly used in the literature, for example, in  Ding et al. (\cite{DC2020}, 2020) and He and Chen (\cite{HC2021}, 2021). We require these conditions in order to apply Lemma \ref{A12}.
\end{remark}
By Theorem \ref{th2}, we can obtain the following accurate asymptotic normality rate for \( f_n(x) - f(x) \).
\begin{theorem}\label{th4}
Suppose the assumptions of Theorem \ref{th2} and (H1)-(H3) hold. Assume that  $2^k \to \infty$  and  $2^k/n \to 0$. Then for any $x\in [0, 1]$, 
$$
    \sup_{-\infty < y < \infty} \left| P\left(\frac{f_n(x) - f(x)}{\sqrt{Var[f_n(x)]}}  \leq y \right) - \Phi(y) \right| = O(n^{-1/9}) \,+\, O(n^{-\gamma-1/2})\, +\, O(\tau_k/\sqrt{n}), \quad \text{as} \quad n\rightarrow \infty,
$$
where $\gamma >0$ and
\begin{equation}\label{eqqq1}
\tau_k =
\begin{cases}
(1/2^k)^{\nu-1/2} & \text{if } 1/2 < \nu < 3/2, \\
\sqrt{k}/2^k & \text{if } \nu = 3/2, \\
1/2^k & \text{if } \nu > 3/2.
\end{cases}
\end{equation}
\end{theorem}
\begin{remark}
Note that by choosing \( \nu > 3/2 \), selecting \( 2^k = O(n^{1/3}) \), and for \( \gamma > 0 \), we obtain the same convergence rate as in the previous theorems, namely
$$
    \sup_{-\infty < y < \infty} \left| P\left(\frac{f_n(x) - f(x)}{\sqrt{Var[f_n(x)]}}  \leq y \right) - \Phi(y) \right| = O(n^{-1/9}), \quad \text{as} \quad n\rightarrow \infty.
$$
\end{remark}
\section{Simulation study}
In this section, we conduct a simulation to study the numerical performance of our main results established in the paper. To achieve this, we use a real case to generate a LENQD sequence from the first-order moving average (MA(1)) process \( X_t = W_t - b W_{t-1} \), a result of an extensive study of annual temperature values measured in Basel from 1755 to 1957 (where \( W_t \) is a Gaussian white noise process with mean \( \mu_W \) and standard deviation \( \sigma_W \)). This time series has been studied by Gilgen (\cite{G2006}, 2006),  more recently by Kheyri et al. (\cite{K2019}, 2019) and most recently by El Alem et al. (\cite{E2025}, 2025). The study showed that the best MA(1) model is the one with \( b = 0.9 \), \( \mu_W = 0 \), and \( \sigma_W = 0.7 \) according to the Akaike Information Criterion (AIC). For the normality test, the Jarque-Bera test was used, yielding a corresponding p-value of 0.839. To test for autocorrelation, the usual test employed in this approach is the portmanteau test. The null hypothesis of no autocorrelation is not rejected, as the p-values exceed the 5\% significance level. Since the sequence \(\{ W_n, n\geq 0 \} \) consists of i.i.d. random variables with standard normal distributions \( N(0, \sigma_W^2) \), and \( X_n = W_n - b W_{n-1} \), it follows that \( \{X_n, n\geq 1 \} \) is a Gaussian sequence with the same distribution \( N(0, (1 + b^2)\sigma_W^2) \). The covariance structure is given by \( \text{Cov}(X_i, X_j) = -b\sigma_W^2 \) if \( |i - j| = 1 \), and \( \text{Cov}(X_i, X_j) = 0 \) otherwise. As a negatively correlated Gaussian sequence \( \{X_n, n \geq 1\} \), it is a NA sequence and, consequently, LENQD. In fact, the finite-dimensional distributions of this sequence are \( (X_1, \dots, X_n) \sim \mathcal{N}_n(0, \Sigma) \), where \( \Sigma = [a_{ij}] \) is an \( n \times n \) covariance matrix with diagonal entries \( a_{ii} = (1 + b^2)\sigma_W^2 \), and off-diagonal entries \( a_{ij} = -b\sigma_W^2 \) for \( i = j-1 \) or \( j+1 \), and \( a_{ij} = 0 \) otherwise. The simulation procedure for the sequence \( \{X_n, n \geq 1\} \) is straightforward.\\

To illustrate how the quality of  the CLT result established in Theorem 2.1 is influenced by the sample size \( n \), we performed \( M = 1000 \) replicates with \( n \in \{100, 300, 500\} \). Next, for these different sample sizes, the Quantile-Quantile  (Q-Q) plot of sample data versus the standard normal distribution is generated to highlight the favorable behavior of our results. This numerical study shows that the goodness of fit improves with an increasing sample size $n$, which is clearly apparent in the Figures \ref{fig1}, \ref{fig2} and \ref{fig3}.

\begin{figure}[h!]\label{fig1}
\begin{center}
\includegraphics[height=5.5cm, width=0.8\textwidth]{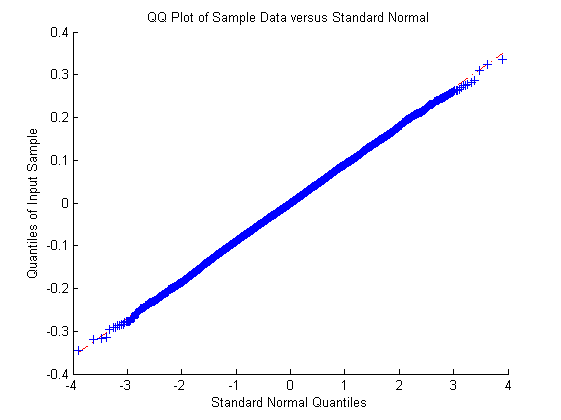}
  \caption{QQ plot of sample data versus the standard normal  for $n=100$}\label{fig1}

 \end{center}
\end{figure}
\begin{figure}[h!]\label{fig2}
\begin{center}
\includegraphics[height=5.5cm, width=0.8\textwidth]{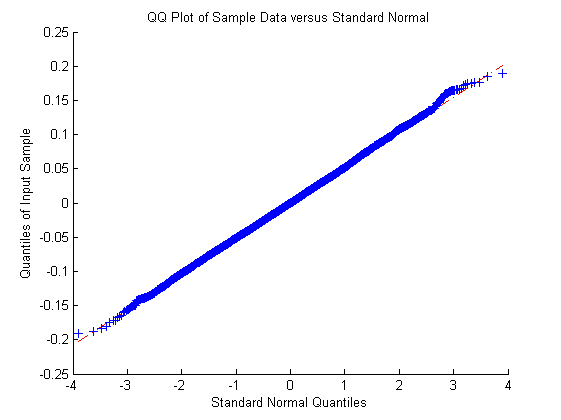}
\caption{QQ plot of sample data versus the standard normal  for $n=300$}\label{fig2}

 \end{center}
\end{figure}
\begin{figure}[h!]
\begin{center}
\includegraphics[height=5.5cm, width=0.8\textwidth]{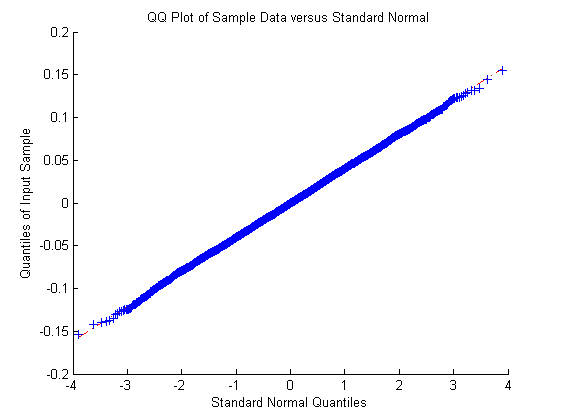}
 \caption{QQ plot of sample data versus the standard normal  for $n=500$}\label{fig3}
 \end{center}
\end{figure}
To show the performance of the wavelet estimator, we compute the empirical distribution function of \( S:=\frac{f_n(x) - E[f_n(x)]}{\sqrt{Var[f_n(x)]}} \) to estimate \( F(y)= P\left(S  \leq y \right) \) and then estimate the maximum value of \( | F(y) - \Phi(y) | =: Error \) for \( y \in [-3, 3] \).  We simulate data from the following model:
\[
Y_i = f(x_i) + \varepsilon_{i}, \quad 1 \leq i \leq n,
\]
where 
\[
f_1(x) = 2x-1, \quad f_2(x) = \sin(2\pi x),  \quad f_3(x) = \exp(-2x), \quad x_i = \frac{i}{n}.
\]
For generating the error sequence \( \{ \varepsilon_{i}; i=1,...,n\} \), we first generate \( (n + 1) \) i.i.d. random variables \( W_i \) drawn from \( N(0, \sigma_W^2) \). Then, we generate the LENQD sequence \( \{\varepsilon_{i}; i = 1, \cdots, n\} \) by \( \varepsilon_{i} = W_i - b W_{i-1} \), where \( b = 0.9 \) and \( \sigma_W = 0.7 \). We choose the scale function \( \phi(x) = I(0 \leq x \leq 1) \), and take \( 2^k = n^{1/3} \) for a given sample size n.  From the above model, we generate the observed data with sample sizes \( n \in \{100, 300, 500\} \), respectively. For each case, the simulation is repeated 1000 times.

\begin{table}[h]\label{table1}
\begin{center}
\caption{The uniform Berry-Esseen bounds for wavelet estimator}\label{table1}
\begin{tabular}{cccc}
\hline
n & 100 & 300 & 500 \\ \hline
\begin{tabular}{c}
$f(x)$  \\ \hline
$2x-1$ \\ 
$\sin(2\pi x)$ \\ 
$\exp(-2x)$
\end{tabular}
& 
\begin{tabular}{l}
 Error \\ \hline
0.0256 \\
0.0149 \\ 
0.0315
\end{tabular}
& 
\begin{tabular}{l}
 Error \\ \hline
0.0232 \\
0.0141\\ 
0.0305
\end{tabular}
& 
\begin{tabular}{c}
 Error \\ \hline
$0.0212$ \\
$0.0133$\\ 
$0.0267$
\end{tabular}
\\ \hline
\end{tabular}
\end{center}
\end{table}
From the results in Table 1, we observe a good fit of our main results on the Berry-Esseen bound of the wavelet estimator for the nonparametric regression model established in the article as \( n \) increases.
\section{Some useful lemmas}
Before giving the proofs of our main results, we need to state some preliminaries lemmas.
\begin{lemma}[Chang and Rao (\cite{CR}, 1989)]\label{A3}
Let $X$ and $Y$ be random variable. For any $a> 0$ we have 
$$
\sup_{-\infty < u < \infty}\left\vert P(X+Y\leq u)-\Phi(u)\right\vert \leq\sup_{-\infty < u < \infty}|P(X\leq u)-\Phi(u)|+\dfrac{a}{\sqrt{2\pi}}+P(|Y|>a).
$$
\end{lemma}
\begin{lemma}[Li et al. (\cite{Li2023}, 2023)]\label{A1}
Let \( p \geq 2 \) and \(\{X_n, n \geq 1\}\) be a sequence of zero mean LENQD rv's with dominating constant $M\geq 1$ and \( E|X_n|^p < \infty \) for each \( n \geq 1 \). Then there exists a positive constant \( C_p \) depending only on \( p \) such that
\[
E \left\vert \sum_{i=1}^n X_i \right\vert^p \leq C_p \left\{ \sum_{i=1}^n E|X_i|^p + M \left( \sum_{i=n_1}^n E X_i^2 \right)^{p/2} \right\}.
\]
\end{lemma}
\begin{lemma}[Li et al. (\cite{Li2023}, 2023)]\label{ChIneq}
Let $X_1, X_2, \cdots, X_{k_n} $ be a sequence of LENQD rv's with dominanting constant $M\geq 1$, and with finite second moments. Then, for any $t\in \mathbb{R}$,
	$$
	\left \vert E\exp\left\{it\sum_{m=1}^{k_n} X_{m}  \right\}- \prod_{m=1}^{k_n} E\exp \left\{itX_{m}  \right\}\right \vert \leq  2t^2\sum_{1\leq i< j\leq k_n}\left \vert Cov\left(X_{i}, X_{j} \right) \right \vert.
	$$
\end{lemma}
\begin{lemma}[Li et al. (\cite{Li2023}, 2023)]\label{A11}
Let \(\{X_n, n \geq 1\}\) be a sequence of LENQD rv's with dominating coefficient $M\geq 1$. If \(\{ \upsilon_n, n \geq 1 \}\) is a sequence of all nondecreasing (or all nonincreasing) functions, then \(\{ \upsilon_n(X_n), n \geq 1 \}\) is also a sequence of LENQD rv's with the same dominating coefficient.
\end{lemma}
\begin{lemma}[Antoniadis et al. (\cite{G1994}, 1994)]\label{A12}
Suppose that (H1)-(H3) hold. Then
\[
\left| E[f_n(x)] - f(x)\right| = O(n^{-\gamma}) + O(\tau_k),\quad \text{where}\; \tau_k \; \text{is defined in (\ref{eqqq1}).}
\]
\end{lemma}

\section{Proofs}
\subsection{Proof of Theorem \ref{MainResult}}
\begin{proof}
We employ Bernstein's big-block and small-block procedure. Partition the set
$\{1, 2, \ldots, n\}$ into $2k_n + 1$ subsets with large blocks of size $p_n$ and small blocks of size $q_n$. Define
\[
p_n = \left[ n^{2/3} \right], \quad q_n = \left[ n^{1/3} \right], \quad k=k_n := \left[ \frac{n}{p_n + q_n} \right] = \left[ n^{1/3} \right], \quad \text{and}\quad Z_{n,i} = \frac{X_i}{V_n}.
\]
For $j=0, 1, \cdots, k-1$, let \( \eta_j \), \( \xi_j \), and \( \zeta_j \) be defined as follows: 
\[ \eta_j := \sum_{i=j(p_n + q_n) + 1}^{j(p_n + q_n) + p_n} Z_{n,i}, \quad
 \xi_j := \sum_{i=j(p_n + q_n) + p_n + 1}^{(j + 1)(p_n + q_n)} Z_{n,i},  \quad \text{and} \quad
 \zeta_k := \sum_{i=k(p_n + q_n) + 1}^n Z_{n,i}. \]
Writing
\[
\Sigma_n := \frac{S_n}{V_n} = \sum_{j=0}^{k-1} \eta_j + \sum_{j=0}^{k-1} \xi_j + \zeta_k =: \Sigma^{'}_n + \Sigma^{''}_n + \Sigma^{'''}_n.
\]
By using Lemma \ref{A3} with $a=2n^{-1/9}$, we can see that
\begin{align}
\sup_{-\infty < y < \infty} \left| P(\Sigma_n \leq y) - \Phi(y) \right| & \leq \sup_{-\infty < y < \infty} \left| P(\Sigma^{'}_n \leq y) - \Phi(y) \right| \notag \\
& + \frac{2n^{-1/9}}{\sqrt{2\pi}} + P(|\Sigma^{''}_n| > n^{-1/9}) + P(|\Sigma^{'''}_n| > n^{-1/9}). \notag
\end{align}

Firstly, we estimate \( E( \Sigma^{''}_n)^2 \) and \( E(\Sigma^{'''}_n)^2 \), which will be used to estimate \( P(\vert \Sigma^{''}_n \vert > n^{-1/9}) \) and \( P(\vert \Sigma^{'''}_n \vert > n^{-1/9}) \). It  follows by the conditions \( \vert X_i \vert \leq c \) and (C2) that  
\begin{equation}\label{eq1}
|Z_{n,i}|\leq \frac{\vert X_i \vert}{V_n}  \leq \frac{c_1}{\sqrt{n}}.
\end{equation}
Combining the definition of LENQD rv's  with the definition of \( \xi_j \), \( j = 0, 1, \ldots, k - 1 \), we can easily prove that \( \{\xi_0, \xi_1, \ldots, \xi_{k-1}\} \) is LENQD. Therefore, it follows from (\ref{eq1}), \( E[Z_{n,i}] = 0 \)  and Lemma \ref{A1} that
$$
E( \Sigma^{''}_n)^2 \leq c_1 \sum_{j=0}^{k-1} E(\xi_j^2) \leq c_2 \frac{k_nq_n}{n} \leq c_3 \frac{q_n}{p_n} = O(n^{-1/3}),
$$
and
\begin{align}\label{eqqq}
E(\Sigma^{'''}_n)^2 &\leq  \frac{c_4}{n}E \left\vert \sum_{i=k(p_n+q_n)+1}^{n} X_i \right\vert^2
\leq \frac{c_5}{n} \sum_{i=k(p_n+q_n)+1}^{n} E(X_i^2)  \notag \\
&\leq \frac{c_6}{n}(n - k_n (p_n + q_n)) \leq c_7 \frac{p_n + q_n}{n} = O(n^{-1/3}).\notag
\end{align}
Consequently, by Markov's inequality, we obtain
\[
P(\vert \Sigma^{''}_n \vert > n^{-1/9}) \leq n^{2/9} E( \Sigma^{''}_n)^2 = O(n^{-1/9}),
\]
\[
P(\vert \Sigma^{'''}_n \vert > n^{-1/9}) \leq n^{2/9} E( \Sigma^{'''}_n)^2 = O(n^{-1/9}).
\]
In the following, we will estimate
\(
\sup_{-\infty < y < \infty} \left| P(\Sigma_n \leq y) - \Phi(y) \right|.
\)
Define
\[
s_n^2 := \sum_{j=0}^{k-1} \text{Var}(\eta_j), \quad \text{and} \quad \Gamma_n := \sum_{0 \leq i < j \leq k-1} \text{Cov}(\eta_i, \eta_j).
\]
Here, we first estimate the growth rate \( |s_n^2 - 1| \). Since \( E[\Sigma_n^2] = 1 \), then
\[
E(\Sigma_n')^2 = E \left[ (\Sigma_n - (\Sigma_n'' + \Sigma_n'''))^2 \right] = 1 + E \left[ (\Sigma_n'' + \Sigma_n''')^2 \right] - 2E \left[ \Sigma_n (\Sigma_n'' + \Sigma_n''') \right].
\]
Thus, it follows that
\begin{align}\label{eqqq}
\vert E(\Sigma_n')^2 - 1\vert &=\left| E \left[ (\Sigma_n'' + \Sigma_n''')^2 \right] - 2E \left[ \Sigma_n (\Sigma_n'' + \Sigma_n''')^2 \right] \right| \notag \\
&\leq E(\Sigma_n'')^2 + E(\Sigma_n''')^2 + 2 \sqrt{E[(\Sigma_n'')^2]E[(\Sigma_n''')^2]}\notag \\ 
&\qquad + 2  \sqrt{E[(\Sigma_n)^2]E[(\Sigma_n'')^2]}  + 2  \sqrt{E[(\Sigma_n)^2]E[(\Sigma_n''')^2]} 
 \notag \\
&=O(n^{-\frac{1}{3}}) + O(n^{-\frac{1}{6}}) = O(n^{-\frac{1}{6}}). \notag
\end{align}
Notice that
\[
s_n^2 = E(\Sigma_n')^2 - 2 \Gamma_n. 
\]
With \( l_j = j(p_n + q_n) \),
\[
2\Gamma_n = 2 \sum_{0 \leq i < j \leq k-1} \sum_{\lambda_1=1}^{p_n} \sum_{\lambda_2=1}^{p_n} \text{Cov}(Z_{n,l_i + \lambda_1}, Z_{n,l_j + \lambda_2}),
\]
but since \( i \ne j \), \( |l_i - l_j + l_1 - l_2| \geq q_n \), it follows that
\begin{align}\label{eqqq}
2|\Gamma_n|& \leq 2 \sum_{1 \leq i < j \leq n, \, j - i \geq q_n} |\text{Cov}(Z_{n,i}, Z_{n,j})|
 \notag \\
&\leq \frac{c_1}{n} \sum_{1 \leq i < j \leq n, \, j - i \geq q_n} |\text{Cov}(X_i, X_j)|
\notag \\
&\leq c_2 \sum_{r = q_n}^{\infty} |\text{Cov}(X_1, X_r)| = O(n^{-\delta/3}) = O(n^{-1/3}). \notag
\end{align}
From the conditions (C1)-(C2),  we can get that
\begin{align}
|s_n^2 - 1|& \leq \vert E(\Sigma_n')^2 - 1\vert +2|\Gamma_n|  \notag \\
& =  O(n^{-1/6}) + O(n^{-1/3}) \notag \\
&= O(n^{-1/6}). \notag
\end{align}
For \( j = 0, 1, \ldots, k - 1 \), let \(\hat{\eta}_j \) be independent rv's and have the same distribution as \( \eta_j \). Define \( G_n = \sum_{j=0}^{k-1} \eta_j \). It can be found that
\begin{align}
\sup_{-\infty < y < \infty} \left| P(\Sigma_n' \leq y) - \Phi(y) \right|
 &\leq \sup_{-\infty < y < \infty} \left| P(\Sigma_n' \leq y) - P(G_n \leq y) \right|  \notag \\
& + \sup_{-\infty < y < \infty} \left| \Phi\left(y/s_n\right) - \Phi(y) \right| + \sup_{-\infty < y < \infty} \left| P(G_n \leq y) - \Phi\left(y/s_n\right) \right|  \notag \\
&:= R_1 + R_2 + R_3.\notag
\end{align}
Let $\varphi(y)$ and $\psi(y)$ be the characteristic functions of $\Sigma_n'$ and $G_n$, respectively. By Esseen inequality [(\cite{Petrov}, 1995), Theorem 5.3], for any $\tau >0$,
\begin{align}
R_1 & \leq \int_{-\tau}^{\tau} \left|\frac{\varphi(y) -  \psi(y)}{y} \right| \, dy + \tau \sup_{-\infty < y < \infty} \int_{|u| \leq c/\tau} \left| P(G_n \leq u + y) - P(G_n \leq y) \right| \, du \notag \\
&:= R_{1n} + R_{2n}. \notag
\end{align}
With \( l_j = j(p_n + q_n) \) and  Lemma \ref{ChIneq}, we have that
\begin{align}
|\varphi(t) - \psi(t)| & = \left| E \left[ \exp \left( it \sum_{j=0}^{k-1} \eta_j \right) \right] - \prod_{j=0}^{k-1} E \left[ \exp(it \eta_j) \right] \right| \notag \\
& \leq ct^2 \sum_{0 \leq i < j \leq k-1} \sum_{\lambda_1=1}^{q_n} \sum_{\lambda_2=1}^{q_n} \left| \text{Cov}(Z_{n,l_i + \lambda_1}, Z_{n,l_j + \lambda_2}) \right| \notag \\
& \leq c_1 t^2 \sum_{1 \leq i < j \leq n, \, j - i \geq q_n} \left| \text{Cov}(X_i, X_j) \right|
 \notag \\
& \leq c_2 t^2 \sum_{r = q_n}^{\infty} \left| \text{Cov}(X_1, X_r) \right| \notag \\
&\leq c_3 t^2 n^{-\delta/3}. \notag
\end{align}
Set \( \tau = n^{ (\delta - 1)/6} \) for \( \delta \geq 1 \). So, we have
\[
R_{1n} = \int_{-\tau}^{\tau} \left|\frac{\varphi(y) -  \psi(y)}{y} \right| \, dy \leq c n^{-\delta/3} \cdot \tau^2 = O(n^{-1/3}).
\]
It follows from the Berry-Esseen inequality (\cite{Petrov}, Theorem 5.7) that
\[
\sup_{-\infty < y < \infty} \left| P \left( G_n/s_n \leq y \right) - \Phi(y) \right| \leq \frac{c}{s_n^3} \sum_{j=0}^{k-1} E \left|\hat{\eta}_j \right|^3 = \frac{c}{s_n^3} \sum_{j=0}^{k-1} E \left| \eta_j \right|^3.
\]
By condition $\vert X_i \vert \leq c$ and Lemma \ref{A1},  we have
\begin{align}
\sum_{j=0}^{k-1} E \left| \eta_j \right|^3 & = \sum_{j=0}^{k-1} E \left[ \left( \sum_{i= j(p_n + q_n)+1}^{j(p_n + q_n)+p_n} Z_{n,i} \right)^3 \right] \notag \\
& \leq \frac{c}{n^{3/2}} \sum_{j=0}^{k-1} E \left[ \left( \sum_{i= j(p_n + q_n)+1}^{j(p_n + q_n)+p_n}  X_{i} \right)^3 \right] \notag \\
& \leq \frac{c_1}{n^{3/2}} \sum_{j=0}^{k-1} \left\{  \sum_{i= j(p_n + q_n)+1}^{j(p_n + q_n)+p_n}  E \left\vert X_{i} \right\vert^3 +M\left(\sum_{i= j(p_n + q_n)+1}^{j(p_n + q_n)+p_n}  E X_{i}^2  \right)^{3/2} \right\} \notag \\
& \leq \frac{c_2}{n^{3/2}} \sum_{j=0}^{k-1} \left( p_n + p_n^{3/2} \right) \leq\frac{c_3kp_n^{3/2}}{n^{3/2}} = O(n^{-1/6}). \notag
\end{align}
Hence, yield that
\begin{equation}\label{R2}
\sup_{-\infty < y < \infty} \left| P \left( G_n/s_n \leq y \right) - \Phi(y) \right|= O(n^{-1/6}).
\end{equation}
Now, since \( s_n \to 1 \) as \( n \to \infty \), it suffices to show that
\begin{align}
\sup_{-\infty < y < \infty} \left| P(G_n \leq u + y) - P(G_n \leq y) \right| & \leq \sup_{-\infty < y < \infty} \left| P \left( G_n/s_n \leq (u + y)/s_n \right) - \Phi \left( (u + y)/s_n \right) \right|
 \notag \\
& \qquad + \sup_{-\infty < y < \infty} \left| P \left( G_n/s_n \leq y/s_n \right) - \Phi \left( y/s_n \right) \right| \notag \\
&\qquad + \sup_{-\infty < y < \infty} \left| \Phi \left( (u + y)/s_n \right) - \Phi \left( y/s_n \right) \right| \notag \\
& \leq 2 \sup_{-\infty < y < \infty} \left| P \left(G_n/s_n \leq y \right) - \Phi(y) \right| \notag \\
&\qquad + \sup_{-\infty < y < \infty} \left| \Phi \left( (u + y)/s_n \right) - \Phi \left( y/s_n \right) \right| \notag \\
& = O(n^{-1/6}) + O \left( \left| u\right|/s_n  \right), \notag
\end{align}
which yields that
\begin{align}
R_{2n} & = \tau \sup_{-\infty < y < \infty} \int_{|u| \leq c/\tau} \left| P(G_n \leq u + y) - P(G_n \leq y) \right| \, du \notag \\
& \leq  c_1/ n^{1/6} + c_2/\tau  \notag \\
& = O(n^{-1/6}) + O(n^{-1/3}) = O(n^{-1/6}), \quad \text{where}\quad \tau = n^{1/3}. \notag
\end{align}
From the inequality (\cite{Petrov}, Lemma 5.2)
\[
\sup_{-\infty < x < \infty} \left| \Phi(\alpha x) - \Phi(x) \right| \leq
 \frac{(\alpha - 1) I(\alpha \geq 1)}{\sqrt{2 \pi e}} +  \frac{(1/\alpha - 1) I(0 < \alpha < 1)}{\sqrt{2 \pi e}},
\]
it can be conclued that,
\begin{align}
R_2 & = \sup_{-\infty < y < \infty} \left| \Phi \left( y/s_n \right) - \Phi(y) \right| \notag \\
& \leq \frac{(s_n - 1) I(s_n > 1)}{\sqrt{2 \pi e}}  + \frac{(1/s_n - 1) I(0 < s_n \leq 1)}{\sqrt{2 \pi e}} 
 \notag \\
&\leq \frac{\max \left( |s_n - 1|, (|s_n - 1|)/s_n \right)}{\sqrt{2 \pi e}} 
 \notag \\
&\leq c_1 \max \left( |s_n - 1|, (|s_n - 1|)/s_n \right) \cdot (s_n + 1) \quad \text{(note that } s_n \to 1\text{)} \notag \\
& \leq c_2 |s_n^2 - 1| = O(n^{-1/6}), \notag
\end{align}
and by (\ref{R2})
$$
R_3= \sup_{-\infty < y < \infty} \left| P(G_n \leq y) - \Phi\left(y/s_n\right) \right| = O(n^{-1/6}).
$$
Therefore, it follows that 
$$
\sup_{-\infty < y < \infty} \left| P(\Sigma_n' \leq y) - \Phi(y) \right| =  O(n^{-1/6}).
$$
Finally, the proof is completed.
\end{proof}
\subsection{Proof of Theorem \ref{th2}}
\begin{proof}
Firstly, under the assumption that \( \int_{\Gamma_{i}} E_k(x, s) \, ds > 0 \) for all \( 1 \leq i \leq n \), we recall from Lemma~\ref{A11} that if \( \{X_n, n \geq 1\} \) is a sequence of LENQD rv's, then \( \{aX_n, n \geq 1\} \) is also a sequence of LENQD rv's for any \( a \in \mathbb{R}^* \). Therefore, for a fixed \( x \), the sequence \( \left\{ \int_{\Gamma_{i}} E_k(x, s) \, ds \, \varepsilon_{i}, \; 1 \leq i \leq n \right\} \) is LENQD with a dominating constant \( M \), as it is a non-decreasing transformation of LENQD rv's.\\
Since $E [\varepsilon_{i}] = 0$, in order to prove Theorem \ref{th2}, it suffices to show that,
$$
f_n(x) - E[f_n(x)]=\sum_{i=1}^n\int_{\Gamma_{i}} E_k(x, s) \, ds\, \varepsilon_{i},
$$
and  
$$ 
Var[f_n(x)]=Var\left[\sum_{i=1}^n\int_{\Gamma_{i}} E_k(x, s) \, ds\, \varepsilon_{i} \right].
$$
It can be checked finally that all the conditions of Theorem \ref{MainResult} are satisfied with $X_i= \int_{\Gamma_{i}} E_k(x, s) \, ds\, \varepsilon_{i}$. This completes the proof of the theorem.
\end{proof}
\subsection{Proof of Theorem \ref{th3}}
\begin{proof}
By applying (\ref{E}) in conjunction with Lemma \ref{A3} for  $a = \frac{\left| \sum_{i=1}^n \left( \int_{\Gamma_i} E_k(x, s) \, ds \right) \varepsilon_i \right|}{\sqrt{Var(f_n(x))}} $, we obtain the following result 
\begin{align}
\sup_{y \in \mathbb{R}} \left| P \left(\frac{f_n(x)-E[f_n(x)]}{\sqrt{Var(f_n(x))}} \leq y \right) - \Phi(y) \right| & \leq \sup_{y \in \mathbb{R}} \left| P \left( \frac{ \sum_{i=1}^n \left( \int_{\Gamma_i} E_k(x, s) \, ds \right)^+ \varepsilon_i}{\sqrt{Var(f_n(x))}} \leq y \right) - \Phi(y) \right|  \notag \\
&\qquad \quad +\;  \frac{\left| \sum_{i=1}^n \left( \int_{\Gamma_i} E_k(x, s) \, ds \right) \varepsilon_i \right|}{\sqrt{2\pi}\sqrt{Var(f_n(x))}}  \notag  \\
& =: I_1 + I_2. \notag
\end{align}
Here we use the fact that $$P \left(  \frac{\left| \sum_{i=1}^n \left( \int_{\Gamma_i} E_k(x, s) \, ds \right)^- \varepsilon_i \right|}{\sqrt{Var(f_n(x))}}  > a \right) = P(\varnothing)=0\quad \text{for the selected}\quad a.$$
Next, observe that by the triangle inequality and the well-known inequality 
\begin{equation}\label{In}
\sup_{x \in \mathbb{R}} |\Phi(ax) - \Phi(x)| \leq \frac{1}{e \sqrt{2\pi}} \left( |a - 1| + \left| a^{-1} - 1 \right| \right),
\end{equation}
we have:
\begin{align}
I_1 & = \sup_{y\in\mathbb{R}} \left| {P} \left( \frac{ \sum_{i=1}^n \left( \int_{\Gamma_i} E_k(x, s) \, ds \right)^+ \varepsilon_i }{\sigma_n^+(x) } \leq \frac{\sigma_n(x)}{\sigma_n^+(x) } y \right) - \Phi(y) \right| \notag \\
&\leq \sup_{y\in\mathbb{R}} \left| {P} \left( \frac{ \sum_{i=1}^n \left( \int_{\Gamma_i} E_k(x, s) \, ds \right)^+ \varepsilon_i }{\sigma_n^+(x) } \leq \frac{\sigma_n(x)}{\sigma_n^+(x) } y \right) - \Phi\left(\frac{\sigma_n(x)}{\sigma_n^+(x) }y\right) \right| + \sup_{y\in\mathbb{R}} \left|\Phi\left(\frac{\sigma_n(x)}{\sigma_n^+(x) }y\right) - \Phi(y)  \right|  \notag 
\end{align}
It can be verified that all the conditions of Theorem \ref{MainResult} are satisfied with
\(
X_i = \left( \int_{\Gamma_i} E_k(x, s) \, ds \right)^+ \varepsilon_i.
\)
Therefore, by applying Theorem \ref{MainResult}, together with inequality (\ref{In}) and Assumption (\ref{A}), it follows that
$$
I_1=O(n^{-1/9})+O(n^{-\kappa})=O(n^{-1/9}).
$$
On the other hand, it is easily seen that 
$$
I_2 \leq \frac{c}{\sqrt{{Var}(f_n(x))}} = O(n^{-1/2}),
$$
because, by the assumption in Theorem \ref{th2}, \( \sqrt{\text{Var}[f_n(x)]} \geq \sigma \sqrt{n}/2 \) for all sufficiently large \( n \). From the statements above, the proof is completed
\end{proof}
\subsection{Proof of Theorem \ref{th4}}
\begin{proof}
By Lemma \ref{A3}, for any $a>0$, we have:
\begin{align}
\sup_{y \in \mathbb{R}} \left| P \left(\frac{f_n(x)-f(x)}{\sqrt{Var(f_n(x))}} \leq y \right) - \Phi(y) \right| & = \sup_{y \in \mathbb{R}} \left| P \left( \frac{f_n(x)- E[f_n(x)]}{\sqrt{Var(f_n(x))}} + \frac{E[f_n(x)]-f(x)}{\sqrt{Var(f_n(x))}} \leq y \right) - \Phi(y) \right| \notag \\
& \leq \sup_{y \in \mathbb{R}} \left| P \left( \frac{f_n(x)-E[f_n(x)]}{\sqrt{Var(f_n(x))}} \leq y \right) - \Phi(y) \right| + \frac{a}{\sqrt{2\pi}}  \notag \\
&\qquad \quad +\; P \left(  \frac{\left|E[f_n(x)]-f(x)\right|}{\sqrt{Var(f_n(x))}}  > a \right)
 \notag \\
&= O(n^{-1/9}) + \frac{a}{\sqrt{2\pi}} + P \left(  \frac{\left|E[f_n(x)]-f(x)\right|}{\sqrt{Var(f_n(x))}}  > a \right).  \notag
\end{align}
Now, it follows from the assumption in Theorem \ref{th2} that \( \sqrt{\text{Var}[f_n(x)]} \geq \sigma \sqrt{n}/2 \) for all sufficiently large \( n \). Then, by Lemma \ref{A12}, we have
$$
 \frac{\left|E[f_n(x)]-f(x)\right|}{\sqrt{Var(f_n(x))}}  \leq \frac{\left|E[f_n(x)]-f(x)\right|}{\sigma \sqrt{n}/2} = O(n^{-\gamma-1/2}) + O(\tau_k/\sqrt{n}).
$$
Hence, there exists a constant $c$ sufficiently large such that $$\frac{\left|E[f_n(x)]-f(x)\right|}{\sqrt{Var(f_n(x))}}  < c\cdot\max\{n^{-\gamma-1/2}, \tau_k/\sqrt{n}\}.$$ Let $a = c\cdot\max\{n^{-\gamma-1/2}, \tau_k/\sqrt{n}\}$, then $P \left(  \frac{\left|E[f_n(x)]-f(x)\right|}{\sqrt{Var(f_n(x))}}  > a \right)=P\left(\varnothing \right)=0$. Therefore, the proof of Theorem \ref{th3} holds.
\end{proof}
\section{Conclusion}
Our main contribution in this paper extends some existing asymptotic results stated for i.i.d., strong mixing, and negative dependent structures to LENQD  ones. We provide the Berry-Esseen bound for LENQD rv's and apply it to obtain the rate of normality approximation for the wavelet estimation of the nonparametric regression function based on LENQD errors. A simulation study was conducted to evaluate the performance of our main results based on a real dataset.  The results obtained in this paper generalize some corresponding results found in the literature.



\end{document}